\documentclass[a4paper,11pt]{amsart} 

\usepackage{latexsym}
\usepackage[vcentermath]{youngtab} \newcommand{\yngt}{\youngtabloid}
\usepackage{url}
\usepackage{amsmath}
\usepackage{amsopn,amssymb}
\usepackage{amsthm}
\usepackage{graphics}

\usepackage[margin=0.8in,font=small,labelfont=bf]{caption} 

\newtheorem{theorem}{Theorem}[section]

\newtheorem{lemma}[theorem]{Lemma}
\newtheorem{corollary}[theorem]{Corollary}
\newtheorem{proposition}[theorem]{Proposition}

\newtheorem{conjecture}[theorem]{Conjecture}

\newcommand{\la}{\lambda}

\DeclareMathOperator{\sh}{sh}

\DeclareMathOperator{\sgn}{sgn}

\DeclareMathOperator{\End}{End}

\DeclareMathOperator{\Tr}{Tr}

\newcommand{\ind}{\!\!\uparrow}
\newcommand{\res}{\!\!\downarrow}

\newcommand{\tsp}{\hskip 0.5pt}

\newcounter{indentedromanrmklistcnt}
	{\setcounter{indentedromanrmklistcnt}{0}%
	\begin{list}{(\Alph{indentedromanrmklistcnt})}{%
		\usecounter{indentedromanrmklistcnt}%
		\setlength{\topsep}{3pt}%
		\setlength{\itemsep}{3pt}%
		\setlength{\itemindent}{-0.0in}
		\setlength{\leftmargin}{0.5in}} 
	}%
	{\end{list}}%

\newcounter{rmklistcnt}
	{\setcounter{rmklistcnt}{0}%
	\begin{list}{(\arabic{rmklistcnt})}{%
		\usecounter{rmklistcnt}%
		\setlength{\rightmargin}{0.0in}%
		\setlength{\topsep}{3pt}%
		\setlength{\itemsep}{3pt}%
		\setlength{\itemindent}{0.5in}%
		\setlength{\leftmargin}{0.2in}}%
	}%
	{\end{list}}%

\title[Vertices of Specht modules]{Vertices of Specht modules and blocks of the
symmetric group}

\author{Mark Wildon}

\begin{document}

\begin{abstract}
This paper studies the vertices, in the sense 
defined by J.~A.~Green, of Specht modules for symmetric groups.
The main theorem gives, for each indecomposable non-projective
Specht module, a large subgroup contained in one
of its vertices. A corollary of this theorem is a new way to
determine the defect groups of
symmetric groups. We also use it to find the
Green correspondents of a particular family of simple Specht modules;
as a corollary, we get 
a new proof of the Brauer 
correspondence for blocks of the symmetric group.
The proof of the main theorem
uses the Brauer homomorphism on modules, as developed
by M.~Brou{\'e}, together with combinatorial arguments
using Young tableaux. 
\end{abstract}

\maketitle
\thispagestyle{empty}

\section{Introduction}
In this paper we apply the methods of
local representation theory to the symmetric group.
Our object is twofold: firstly to prove Theorem~\ref{thm:main} below
on the vertices of Specht modules, and secondly to use
this theorem to give short proofs of two earlier results
on the blocks of symmetric groups. Specifically, we 
determine their defect groups, and how blocks of symmetric groups relate
under
the Brauer correspondence to blocks of local subgroups of symmetric groups.
Our main theorem applies to all Specht modules, but is
strongest for partitions
of the form $(n-m,m)$: we  discuss the special case $(n-2,2)$ 
at the end of the paper.

Vertices were first defined in an influential
paper of J.~A.~Green \cite{Green59}. We recall his definition here.
Let~$G$ be a finite group and let $F$ be a field of
prime characteristic $p$. Let~$M$ be an indecomposable
$FG$-module. A subgroup~$Q$ of $G$ is said to be a 
\emph{vertex} of $M$ if there is an indecomposable 
$FQ$-module~$N$
such that $V$ is a summand of the induced module $N\ind_Q^G$,
\emph{and} $Q$ is minimal with this property.
By \cite[page 435]{Green59}, the vertices of $M$ are $p$-groups, and
any two vertices of $M$ are conjugate in $G$.
The module~$N$ is well-defined
up to conjugacy in $N_G(Q)$; it is referred to as the
\emph{source} of~$V$. 

Despite the central role played by vertices
in open problems in modular representation
theory, such as Alperin's Weight Conjecture \cite{AlperinWeights},
little is known about
the vertices of `naturally occurring' modules, such as
Specht modules for symmetric groups. 
See \S 3 below
for the definition of Specht modules, and other prerequisite
results concerning tableaux and blocks of symmetric groups.
We recall here that if~$\lambda$ is
a partition of $n$, then
the Specht module $S^\lambda$, defined over 
a field of zero characteristic, affords the ordinary irreducible
character of the symmetric group~$S_n$ canonically
labelled by~$\lambda$.
When defined over
fields of prime characteristic, Specht modules usually fail
to be  irreducible.
However, by Theorem~\ref{thm:indec},
they are usually indecomposable. Our main result
is a step towards finding their vertices.

\begin{theorem}\label{thm:main}
Let $\lambda$ be a partition of $n$ and let $t$ be a $\lambda$-tableau.
Let $H(t)$ be the subgroup of the
row-stabilising group of $t$ which
permutes, as blocks for its action, the entries of columns of equal
length in $t$.
If the Specht module~$S^\lambda$, defined over a field
of prime characteristic $p$, is indecomposable, then it has a vertex
containing a Sylow $p$-subgroup of $H(t)$.
\end{theorem}

For example,  if $\lambda = (8,4,1)$ and 
\[ 
\raisebox{0.31in}{$t =$}\ \includegraphics{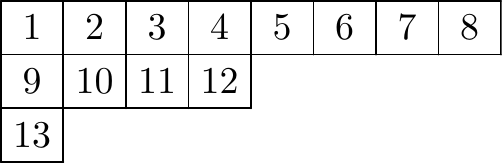} \] 
then the row-stabilising group of $t$ is 
$S_{\{1,2,3,4,5,6,7,8\}} \times S_{ \{9, 10,11,12\} }$ and
$H(t)$ is generated by the permutations
\[ (2,3,4)(10,11,12), (2,3)(10,11), (5,6,7,8), (5,6). \]

Our proof of Theorem~\ref{thm:main}, given in \S 5 below,
uses the Brauer homomorphism
on modules, as developed by M.~Brou{\'e} in~\cite{Broueperm}.
We briefly state the main results we need  from his work in~\S 2. 
We also use a combinatorial result which refines the Standard Basis
Theorem on Specht modules: 
see Proposition~\ref{prop:straightening}. 

We single out the following corollary of Theorem~\ref{thm:main}

\begin{corollary}\label{cor:firstRow}
If the Specht module $S^\lambda$, defined over a field
of characteristic~$p$, is indecomposable,
then it has a vertex containing a Sylow $p$-subgroup of~$S_{\lambda_1 - \lambda_2}$.
\hfill$\qed$
\end{corollary}

In \S 6 and \S 7 we use Corollary~\ref{cor:firstRow} to give new proofs
of two results
on the block theory of the symmetric group. 
We shall suppose in these
sections that the reader has some familiarity with block theory:
see \cite[Chapter~4]{Alperin} for an introduction. 
We recall here 
that if $B$ is a $p$-block of
the finite group $G$ then, when thought of as an $F(G \times G)$-module
with the action $x(g,g') = g^{-1}xg'$ for $x \in B$ and  $g,g' \in G$, 
$B$ has a vertex of the form
\[ \Delta D = \{ (g,g) : g\in D \} \]
for some subgroup $D$ of $G$. We say that $D$ is a \emph{defect group}
of the block~$B$.
By \cite[\S13, Theorem~5]{Alperin}, 
if $M$ is an indecomposable
module lying in a block~$B$
then $M$ has a vertex contained in a defect group of~$B$.

Our results are obtained
by considering a particular family of Specht modules.
Given $w \in \mathbf{N}$ and
a partition $\gamma=(\gamma_1,\ldots,\gamma_k)$ which is a $p$-core,
let
\begin{equation}
\label{eq:gamma} \gamma + wp =  (\gamma_1 + wp, \gamma_2, \ldots, \gamma_k). 
\end{equation}
We shall say that the partitions $\gamma + wp$ are \emph{initial}.
In \S 6 we use Corollary~\ref{cor:firstRow} to 
determine the vertices of Specht modules
labelled by initial partitions. This
gives a new way  to determine the defect groups of blocks
of the symmetric group. 
The ideas in this proof can
also be used to give a short proof of Brauer's Height Zero Conjecture
for the symmetric group. We explain this in~\S 6.1.
 
In  \S 7 we
find the Green correspondents of Specht modules
labelled by initial partitions. 
This leads to a new
way to determine the behaviour of blocks
of the symmetric group under the Brauer correspondence. 
(This was first decided by
M.~Brou{\'e} in \cite{BroueBlocks}.)
By Lemma~\ref{lemma:simple}, each Specht module $S^{\gamma + wp}$ is simple,
so our result is also a first step in finding the Green correspondents
of the simple modules of the symmetric groups. 
For other results on the vertices of particular simple
modules, see \cite{DKZ} and \cite{MZ}.

Most of the existing work on the vertices of Specht modules
has been on Specht modules labelled
by partitions of the form $(n-m,1^m)$. Their vertices
were found by the author in \cite[Theorem~2]{WildonCycSpechtVertices} 
in the case where the
field characteristic does not divide~$n$. The remaining
case was solved in~\cite{MurphyPeel} for fields of characteristic~$2$; 
it is an open problem when $m \ge 2$ for fields of odd characteristic.

We end in \S 8 by using Theorem~\ref{thm:main} to find
the vertices of the Specht modules $S^{(n-2,2)}$ defined over
fields of odd characteristic. The harder case of
characteristic $2$ was recently solved by Danz and Erdmann in
\cite{DanzErdmann}. Theorem~1.1 does, of course, give some
useful information about the vertices of Specht modules labelled
by arbitrary two-part partitions, but we shall 
not attempt
to pursue the problem any further in this paper.

\section{The Brauer homomorphism}
Let $G$ be a 
finite group and let $F$ be a field of prime characteristic $p$.
Let~$M$ be an $FG$-module. 
For $Q \le G$, let 
$M^{Q}$ denote the subspace of~$M$ consisting of those vectors
fixed by every element of $Q$. 
Given subgroups $R \leq Q \leq G$, we define the \emph{relative trace} map
$\Tr_{R}^{Q} : M^{R} \rightarrow M^{Q}$ by
\[ \Tr_{R}^{Q}(x) = \sum_{i=1}^m xg_i \] 
where  $g_1, \ldots, g_m$ is a transversal for the right cosets
of $R$ in $Q$.
The \emph{Brauer quotient}
of $M$ with respect to $Q$ is the quotient space
\[ M(Q) = M^Q \bigl/ \sum_{R < Q} \Tr_R^Q V^R . \]
The \emph{Brauer homomorphism} with respect to $Q$ is the quotient map
$M^Q \twoheadrightarrow M(Q)$. An easy calculation shows that both $M^Q$ and
$\sum_{R < Q} \Tr_R^Q V^R$ are $N_G(Q)$-invariant, and so $M(Q)$
is a module for $FN_G(Q)/Q$. 

The next theorem shows how the
Brauer homomorphism may be used to gather information about vertices. 
It is proved in \cite[(1.3)]{Broueperm}.

\begin{theorem}\label{thm:Brauer}
Let $G$ be a finite group, let $F$ be a field of prime characteristic,
and let $M$ be an indecomposable $FG$-module. Let $Q$ be a subgroup of~$G$.
If $M(Q) \not= 0$ then
$M$ has a vertex containing $Q$.\hfill$\Box$
\end{theorem}

We shall also use
the following theorem, which combines results from Theorem~3.2
of \cite{BroueBlocks} and Exercise~27.4 of \cite{Thevenaz}.

\begin{theorem}\label{thm:Brauerpperm}
Let $G$ be a finite group, let $F$ be a field of prime characteristic,
and let $M$ be an indecomposable $FG$-module with trivial source.

\emph{(i)} If $Q$ is a subgroup of $G$ then $M(Q) \not= 0$ if and only
if $Q$ is contained in a vertex of $M$. 

\emph{(ii)} If $Q$
is a vertex of $M$ then $M(Q)$ is a projective $FN_G(Q)/Q$-module.
Moreover, when
regarded as an $FN_G(Q)$-module, $M(Q)$
is the Green correspondent of $M$.\hfill$\Box$
\end{theorem}

Theorem~\ref{thm:Brauerpperm} cannot be extended to
modules which do not have trivial source.
For example, if $G$ is cyclic of order $4$ and $M$ is the unique
indecomposable $\mathbf{F}_2G$-module of dimension~$3$, then $M(G) = 0$,
even though~$M$ has~$G$ as its vertex.
It is an interesting feature of our proof of Theorem~\ref{thm:main}
that we successfully apply the Brauer homomorphism to modules which
are---in most cases---not trivial source.


\section{Background results on the symmetric group}
In this section we collect the 
prerequisite definitions and results
we need from the representation theory of the symmetric group. 

\subsection{Tableaux}
Let $\lambda$ be a partition of $n$.
A \emph{$\lambda$-tableau} is an assignment of the numbers 
$\{1,2,\ldots,n\}$ to
the boxes of the Young diagram of
 $\lambda$, so that each box has a different entry. 
We say that a
$\lambda$-tableau is \emph{row-standard}
if its rows are increasing when read from left to right, 
and \emph{column-standard} if its columns are increasing
when read from top to bottom. A tableau that
is both row-standard and column-standard is
said to be \emph{standard}.

If $u$ is a tableau, then we denote by $\overline{u}$
the row-standard tableau obtained from $u$ by sorting
its rows in increasing order. We say that $\overline{u}$
is the \emph{row-straightening} of $u$. For example, if
\[ u = \young(4761,253,8) \quad \text{then} \quad \overline{u} =
\young(1467,235,8). \]

Of the many ways to order the set of standard tableaux,
the most fundamental is the dominance order. It will
be useful to define this order on the larger set 
of row-standard tableaux.
First though, we must
define the
dominance order on compositions: if 
$\lambda = (\lambda_1, \ldots, \lambda_\ell)$
and $\mu = (\mu_1, \ldots, \mu_m)$ 
are compositions of $n$, then we say that $\lambda$ \emph{dominates}
$\mu$, and write $\lambda \unrhd \mu$, if
\[ \lambda_1 + \cdots + \lambda_r \ge \mu_1 + \cdots + \mu_r \]
for all $r \in \mathbf{N}$. (If $r$ exceeds the number of parts of
$\lambda$ or $\mu$, then take the corresponding part to be $0$.)
If $t$ is a row-standard tableau, then we denote by $\sh(t^{\le i})$
the composition recording the number of entries $\le i$
in each row of $t$.
For example if 
$t = \overline{u}$ where~$u$ 
is as above,
then $\sh(t^{\le 8}) = (4,3,1)$ and $\sh(t^{\le 5}) = (2,3,0)$.
Finally, if~$\lambda$ is a partition of $n$ and $s$ and $t$
are row-standard $\lambda$-tableaux, then
we say that
$s$ \emph{dominates} $t$ if
\[ \sh\bigl( s^{\le i}\bigr) \unrhd  \sh \bigl( t^{\le i} \bigr) \]
for all $i$ with $1 \le i \le n$. 
Following the usual convention, we shall reuse the~$\unrhd$ symbol for
the dominance order on row-standard tableaux.

\subsection{Specht modules}

We briefly recall the definition of the Specht module~$S^\lambda$
as a submodule of the Young permutation module~$M^\lambda$.
The reader is referred to \cite{James} for examples and 
further details.

Let $\lambda$ be a partition of $n$. Given a $\lambda$-tableau
$t$, we obtain the associated \emph{tabloid} $\mathbf{t}$ by disregarding
the order of the elements within the rows of  $t$. For example, if
\[ t = \young(4761,253,8) \quad \hbox{then} \quad 
\mathbf{t} = \yngt(4761,253,8) = \yngt(1467,235,8) = \ldots \quad \hbox{etc.}\]

\smallskip
\noindent The natural action of $S_n$ on the set of $\lambda$-tableaux
gives rise to a well-defined action of $S_n$ on the set of $\lambda$-tabloids. 
We denote the
 associated permutation representation of $S_n$
 by $M^\lambda$; it is the 
 \emph{Young permutation module} corresponding
to~$\lambda$. If we need to emphasise that the ground ring is~$R$,
then we shall write $M^\lambda_R$.
For example, $M^{(n-1,1)}_{\mathbf{Z}}$ affords the natural 
integral representation
of~$S_n$ as $n \times n$ permutation matrices.

Given a 
$\lambda$-tableau $t$, we let $C(t)$ be the subgroup of $S_n$
consisting of those elements which fix setwise the
columns of $t$. The \emph{polytabloid} corresponding
to $t$ is the element $e_t$ of $M^\lambda$ defined by
\[ e_t = \sum_{g \in C(t)} \textbf{t} g \sgn(g). \]
The \emph{Specht module} $S^\lambda$ is defined to be
the submodule of $M^\lambda$ spanned by the $\lambda$-polytabloids.
Again, we write $S^\lambda_R$ if we need to emphasise the ground ring.
An easy calculation shows that if $h \in S_n$ then
$(e_t)h = e_{th}$, and so $S^\lambda$ is cyclic, generated by
any single polytabloid. 
Moreover, if $h \in C(t)$ then
\begin{equation}
\label{eq:sign} 
(e_t)h = \sgn(h) e_t.
\end{equation}
It follows easily from \eqref{eq:sign}
that $S^\lambda$ is linearly spanned 
by the polytabloids~$e_t$ for which $t$ is column-standard. 
More is true:
in the statement of the following
theorem, if $t$ is a standard tableau, then we say that
$e_t$ is a \emph{standard polytabloid}.

\begin{theorem}[Standard Basis Theorem]\label{thm:std}
The standard $\lambda$-polytabloids form a $\mathbf{Z}$-basis
for the integral Specht module~$S^\lambda_{\mathbf{Z}}$.\hfill$\Box$
\end{theorem}

A short proof
of the Standard Basis Theorem, attributed to J.~A.~Green,
was given in \cite[\S 3]{Peel75}. It is presented with some
simplifications in \cite[Chapter~8]{James}. The corresponding
result for Specht modules defined over fields is an immediate
corollary. 

The next theorem gives two sufficient conditions for
a Specht module to be indecomposable.

\begin{theorem}\label{thm:indec}
Let $F$ be a field of prime characteristic $p$ and
let $\lambda$ be a partition of~$n$. If $p > 2$, 
or if $p=2$ and the parts of $\lambda$ are distinct, then
$S^\lambda_F$ is indecomposable.
\end{theorem}

\begin{proof}
When $p > 2$ it is well known that $\End_{FS_n}(S^\lambda) \cong F$.
(This result has a particularly short proof using the alternative definition
of polytabloids as polynomials: see  \cite[Theorem~4.1]{Peel75} 
or \cite{Specht1935}.
For stronger results in this direction, see \cite[Chapter~13]{James}.)
When the second condition holds, it follows
from Theorems~4.9 and~11.1 in \cite{James} 
that the top of $S^\lambda$ is simple.
Hence in both cases $S^\lambda_F$ is indecomposable.
\end{proof}

When $p=2$, it is 
possible for Specht modules to be decomposable. G.~D.~James
gave the first example in \cite{JamesCounters} where he
showed that $S^{(5,1,1)}_{\mathbf{F}_2}$ is decomposable.
In \cite{MurphyDec}, G.~M.~Murphy showed that this was the first
of infinitely many examples by giving
a necessary
and sufficient condition for the Specht module
$S^{(2m+1-r,1^r)}_{\mathbf{F}_2}$ to be decomposable. In Proposition~3.3.2 of
\cite{WildonThesis},
the author used 
Theorem~2 in \cite{WildonCycSpechtVertices} 
 to give a shorter proof
of Murphy's result.
It is an open
question whether there are any decomposable Specht modules
other than those found by Murphy.

\subsection{Blocks of symmetric groups}
The blocks of symmetric groups are described
by a theorem which seems destined
to remain forever known as Nakayama's Conjecture.
In order to state it we must first
recall some definitions. 

Let $\lambda$ be a partition. 
A \emph{$p$-hook} in $\lambda$ is a connected part of the 
rim of the Young diagram of $\lambda$ consisting of exactly $p$ boxes,
whose removal leaves the diagram of a partition. 
By repeatedly stripping off $p$-hooks from~$\lambda$
we obtain the \emph{$p$-core} of $\lambda$; the number
of hooks we remove is the \emph{weight} of $\lambda$.
For an example, see Figure~1.
Often it is best to perform these operations using
James' \emph{abacus}: for a description of how to use this piece
of apparatus, see \cite[pages 76--78]{JK}.
For instance, it is easy to prove \emph{via} the abacus
that the $p$-core of a partition is well defined,
something which is otherwise
not at all obvious.

\begin{figure}[b]
\includegraphics{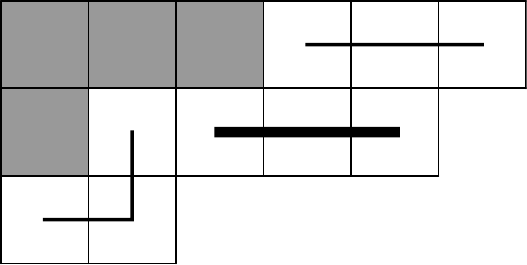}
\caption{The $3$-core of $(6,5,2)$ is $(3,1)$. The thick line
indicates a $3$-hook in $(6,5,2)$; the other two lines
show $3$-hooks of partitions obtained \emph{en route} to
the $3$-core.}
\end{figure}

\vbox{
\begin{theorem}[Nakayama's Conjecture]\label{Thm:Nakayama}
Let $p$ be a prime.
The $p$-blocks of the symmetric group $S_n$ are labelled by pairs $(\gamma, w)$, where
$\gamma$ is a $p$-core and $w \in \mathbf{N}_0$ is the associated weight,
such that $\left|\gamma\right| + wp = n$.
Thus $S^\lambda$ lies
in the block labelled by $(\gamma, w)$ if and only if $\lambda$
has $p$-core $\gamma$ and weight $w$.\hfill$\Box$
\end{theorem}}

Many proofs of Nakayama's Conjecture are now known. A particular
elegant proof was given by Brou{\'e} in \cite{BroueBlocks} using Brauer pairs.
Proposition~2.12 in \cite{BroueBlocks} states the following result
describing the defect groups of blocks of symmetric groups. We shall use
vertices to give
an alternative proof in \S 6 below.

\begin{theorem}\label{thm:SymDefectGroups}
Let $p$ be a prime. If $B$ is a $p$-block of $S_n$ of weight $w$ then the defect
group of~$B$ is a Sylow $p$-subgroup of $S_{wp}$.\hfill$\Box$
\end{theorem}

\medskip

\section{A straightening rule}
The object of this section is to prove the
refinement of the Standard Basis Theorem (Theorem~\ref{thm:std})
stated in  Proposition~\ref{prop:straightening} below.
It seems slightly surprising that this proposition is not
already known; since it appears to be the sharpest 
possible result in  its direction, the author believes that it
well worth putting it on record. 

\medskip

\vbox{
\begin{proposition}\label{prop:straightening} 
Let $\lambda$ be a partition. 
If
$u$ is a column-standard $\lambda$-tableau then its
row-straightening
$\overline{u}$ is a standard tableau. Moreover, 
in $S^\lambda_\mathbf{Z}$,
\[ e_u = e_{\overline{u}} + x \]
where $x$ is an integral linear combination of standard polytabloids~$e_v$ for tableaux 
$v$ such that $\overline{u} \rhd v$.
\end{proposition}}



\begin{proof}
By construction, 
$\overline{u}$ is row-standard, so 
it suffices to show that $\overline{u}$ is also
column-standard.
Let~$\overline{u}_{i,j}$
denote the entry of $\overline{u}$ in row $i$
and column $j$. Suppose, 
for a contradiction, that~\hbox{$\overline{u}_{i+1,j} < \overline{u}_{i,j}$}.
Let
\begin{align*} 
A &= \{ \overline{u}_{i+1,1}, \ldots, \overline{u}_{i+1,j} \}, \\ 
B &= \{ \overline{u}_{i,j}, \ldots, \overline{u}_{i,\lambda_i} \}.
\end{align*}
The entries of~$A$ lie in row $i+1$ of~$u$ and
the entries of~$B$  lie in row~$i$ of~$u$. 
Since
$|A| + |B| = \lambda_i +1$, there exists $a \in A$
and $b\in B$ such that~$a$ and~$b$ appear in the same column
of $u$. 
But then
\[ a \le \overline{u}_{i+1,j} < \overline{u}_{i,j} \le b \]
which contradicts the hypothesis
that $u$ is column-standard.

To prove the second part of the proposition, it will be useful to define the
dominance order on tabloids: given $\lambda$-tabloids
$\mathbf{s}$ and $\mathbf{t}$ corresponding to the tableaux $s$
and $t$ respectively, we set
$\mathbf{s} \tsp \unrhd\tsp \mathbf{t}$ if and only if 
$\overline{s} \tsp \unrhd \tsp \overline{t}$.
%
We shall  need
Lemma~8.3 in \cite{James}, which states that if
$t$ is a column-standard tableau, then
\begin{equation}
\label{eq:dom}
e_t = \mathbf{t} + y
\end{equation}
where $y$ is an integral linear combination of
tabloids $\mathbf{v}$ such that $\mathbf{t} \rhd \mathbf{v}$.

By the Standard Basis Theorem (Theorem~\ref{thm:std})
there exist integers $\alpha_v \in \mathbf{Z}$ such that
\[ e_u = \sum \alpha_v e_v \]
where the sum is over all standard tableaux $v$. 
Let $w$ 
be maximal in the
dominance order such that
$\alpha_w \not =0$.
By \eqref{eq:dom}, applied with $t = w$, 
the tabloid~$\mathbf{w}$
appears with
coefficient 
$\alpha_w$
in $e_u$.  
Another application
of~\eqref{eq:dom}, this time with $t = u$,
shows that $\mathbf{u} \unrhd \mathbf{w}$. 
Hence
\begin{equation}
\label{eq:last}
 e_u = \alpha_{\overline{u}} e_{\overline{u}} + x 
 \end{equation}
where $\alpha_{\overline{u}} \in \{0, 1\}$ and
$x$ is an integral linear combination of standard polytabloids~$e_v$ for tableaux 
$v$ such that $\overline{u} \rhd v$. 

By \eqref{eq:dom},
 the tabloid $\mathbf{u}$ appears in $e_u$ with coefficient $1$.
Again by \eqref{eq:dom}, this tabloid cannot appear in the
summand $x$ in \eqref{eq:last}. It follows that $\alpha_{\overline{u}} = 1$,
as required.
\end{proof}

In addition to our main theorem,
Proposition~\ref{prop:straightening} may also be used to give a short
proof of Lemma~2.1 in G.~E.~Murphy's paper \cite{Murphy}.
It is interesting to note that
the special case of the first part
in which $\lambda$ is a rectangular partition is given---as an
exercise in the Pigeonhole Principle---in \S 10.7 of~\cite{Cameron}.


%



%

%


%

\section{Proof of the main theorem}


%

\smallskip
We are now ready to prove Theorem \ref{thm:main}.
It is clear that we incur no loss in generality if we make 
a specific choice for the $\lambda$-tableau $t$. We
shall take $t$ to be 
the greatest $\lambda$-tableau in the dominance order;
thus if $\lambda = (\lambda_1, \ldots, \lambda_k)$ and
$1 \le j \le k$, then
the entries of $t$ in its $j$th row are
\[ R_j = \{ \la_1+\cdots+\la_{j-1}+1, \ldots, 
\la_1+\cdots+\la_j \}. \]
For example, if $\lambda = (8,4,1)$
then $t$ is the tableau shown after the statement
of Theorem~\ref{thm:main}.
%

We shall show that if $P$ is a Sylow~$p$-subgroup of $H(t)$,
then $S^\lambda(P) \not= 0$.
In the first step, we show that $e_t \in (S^\lambda)^P$.




\vbox{
\begin{lemma}\label{lemma:fixed}
The polytabloid $e_t$ is fixed by every
permutation in $H(t)$.
\end{lemma}

\begin{proof}
Let $h \in H(t)$. By definition, $h$ permutes the columns of $t$
as blocks for its action,  so
$C(t)^h = C(t)$. Moreover, $\textbf{t} h = \textbf{t}$. Hence
\[ e_t h = \sum_{g \in C(t)} \textbf{t}gh \sgn(g) = 
           \sum_{x \in C(t)^h} \textbf{t}hx \sgn(x) = 
           \sum_{x \in C(t)} \textbf{t}x  \sgn(x) = e_t, \]
as required. 
\end{proof}}


The remaining step is given by the case $Q = P$ of the following
lemma. The full strength of this lemma will be
used in \S 7 below. 

\begin{lemma}\label{lemma:kernel}
Let $Q$ be a $p$-subgroup of $H(t)$.
The polytabloid $e_t$ is not contained
in the kernel
of the Brauer homomorphism from $(S^\lambda)^Q$
to $S^\lambda(Q)$.
\end{lemma}

\begin{proof}
Let 
\[ U = \sum_{R < Q} \Tr_R^Q (S^\lambda)^R \]
be the kernel of the Brauer homomorphism with respect to $Q$.
In the sum defining $U$
it suffices to take only those subgroups $R$ which are
maximal subgroups of $Q$, for if $R' < R < Q$ then
\[ \Tr_{R'}^Q (S^\lambda)^{R'} = \Tr_R^Q \Bigl( \Tr_{R'}^R (S^\lambda)^{R'} \Bigr) 
\subseteq
\Tr_R^P (S^\lambda)^R. \]
Hence, if $V$ is the subspace of $S^\lambda$ defined by
\[ V = \left< e_s + e_s g + \cdots + e_s g^{p-1} : 
\text{$s$ a standard $\lambda$-tableau, $g \in Q$} \right>, \]
then  $U$ is contained in $V$. 
We shall use Proposition~\ref{prop:straightening} to show that~$e_t \not\in V$.

Suppose that
there is a standard $\lambda$-tableau $s$,
and a permutation $g \in Q$,
such that $e_t$ appears with a non-zero coefficient in the expression of
\begin{equation}
\label{eq:sum}
e_s + e_s g + \cdots + e_s g^{p-1}
\end{equation}
as a $\mathbf{F}_p$-linear combination of
standard polytabloids. Choose~$i$ such that
$e_t$ appears
with a non-zero coefficient
in the 
expression of $e_sg^i$. 

Let~$u$ be the 
column-standard tableau
whose columns agree setwise with $sg^i$; by \eqref{eq:sign},
$e_u = \pm e_{sg^i}$.
It therefore follows from Proposition~\ref{prop:straightening}
that
\[ 
e_{sg^i} = \pm e_{\overline{u}} + x
\]
where $x$ is a $\mathbf{F}_p$-linear combination of
polytabloids $e_v$ for tableaux $v$ such that $\overline{u} \rhd v$.
Since~$t$ is the greatest tableau in the dominance
order, and the standard polytabloids are linearly independent,
we must have $t = \overline{u}$.

If two elements lying in the same row of $t$
appear in the same column of $sg^i$, then these
elements 
appear in different rows of $\overline{u}$. Hence 
$t \not= \overline{u}$, a contradiction.
We may therefore assume that,
for each $j$, the elements of~$R_j$ appear in 
different columns of $sg^i$. Since $g$ permutes
the elements of each set~$R_j$, it follows that
for each $j$, the elements of $R_j$ appears in different
columns of~$s$.
Since $s$ is standard, row~$j$ of~$s$ must consist
exactly of the elements of~$R_j$; that is, $s = t$.
But, by Lemma~\ref{lemma:fixed}, 
$e_t g = e_t$, and so
\[ e_s + e_s g + \ldots + e_s g^{p-1} = pe_s = 0. \]
This contradicts our assumption that $e_t$ has a non-zero
coefficient in \eqref{eq:sum}. Therefore $e_t \not\in V$, as required.
\end{proof}

Lemmas~\ref{lemma:fixed} and~\ref{lemma:kernel}
imply that $S^\lambda(P)\not=0$. Theorem~\ref{thm:main}
now follows from Theorem~\ref{thm:Brauer}(i).




\smallskip

\section{Defect groups of blocks of the symmetric group}
We now apply Corollary~\ref{cor:firstRow} 
to the initial
partitions defined in~\eqref{eq:gamma} 
to give a new proof of Theorem~\ref{thm:SymDefectGroups} 
on
the defect groups of the symmetric group.
 
Throughout this section, we denote by $[m]_p$ the highest
power of $p$ dividing the natural number $m$.
We shall need the following general
result from block theory, which
connects Brauer's original definition of the 
defect of a block with
our definition \emph{via} vertices. For a proof,
see \cite[Theorem 61.8]{Dornhoff}.

\begin{theorem}\label{thm:dimrestriction}
Let $B$ be a $p$-block of a finite group $G$. Let $p^a$
be the 
highest power of $p$ dividing $|G|$. 
Suppose that 
the defect groups of $B$
have order $p^d$.
If $\chi$ is an irreducible character lying
in $B$, then $p^{a-d}$ divides $\chi(1)$.
Moreover, there is an irreducible character $\chi$ lying in $B$
such that $[\chi(1)]_p = p^{a-d}$. 
\hfill$\Box$
\end{theorem}

We shall also need a companion result, which has an entirely combinatorial
proof.

\begin{lemma}\label{lemma:weightupperbound}
Let $\gamma$ be a $p$-core and let $w\in \mathbf{N}$.
If $[n!]_p = p^a$ and $[(wp)!]_p = p^b$ then
\[ [\dim S^{\gamma + wp}]_p = p^{a-b}. \]
Moreover, if $\mu$ is any other partition with $p$-core $\gamma$
and weight $w$ then 
\[ [\dim S^\mu]_p \ge p^{a-b}. \]
\end{lemma}

\begin{proof}
We shall use the $p$-quotient of a partition, as defined
in \cite[2.7.29]{JK}. 
By \cite[Theorem~2.7.37]{JK}, if $\lambda$ is a partition
with $p$-quotient $(\lambda(0),\ldots, \lambda(p-1))$ then there is a bijection
between
hooks in $\lambda$ of length divisible by $p$ and 
the hooks of the $\lambda(i)$; a hook of length $rp$ in $\mu$
corresponds to a hook of length~$r$ in one of the $\mu(i)$.
The reader may care to verify this for
the partition $(6,5,2)$ shown in Figure~1, which has
$((2),\varnothing,(1))$ as one of its $3$-quotients.

The partition $\gamma + wp$ has $((w), \varnothing, \ldots, \varnothing)$
as a $p$-quotient, and so its $p$-hooks have lengths
$p$, $2p$, \ldots, $wp$. (This can also be seen directly
from its partition diagram.) Hence the highest power of $p$ dividing
the product of the hook-lengths of $\gamma + wp$ is $[(wp)!]_p = p^b$.
The first part of the lemma
now follows
from the Hook Formula for the dimension of Specht modules
(see \cite[Theorem~2.3.21]{JK}).

The second part of the lemma will follow from the Hook Formula if we can show
that the highest power of $p$
dividing the product of the hook-lengths of $\mu$ is at most $p^b$.
Let $(\mu(0), \ldots, \mu(p-1))$ be a $p$-quotient of $\mu$.
Suppose that $\mu(i)$ is a partition of $c_i$.
Writing $h_\alpha$ for the hook-length on a node $\alpha$ of a partition, 
we have
\[ \prod_{\alpha \in \mu} [h_\alpha]_p  = p^w \prod_{i=0}^{p-1}
\prod_{\alpha \in \mu(i)} [h_\alpha]_p \\
 = p^w \prod_{i=0}^{p-1} \left[ \frac{c_i!}{\dim S^{\mu(i)}} \right]_{\!p}.
\]
Rearranging and substituting $[w!]_p$ for $p^{b-w}$ we get
\begin{equation}\label{eq:heightineq}
p^b \,\bigl/\, \prod_{\alpha \in \mu} [h_\alpha]_p = 
{\left[ \binom{w}{c_0,c_1, 
 \ldots, c_{p-1}} \right]_{\!p}}\, \prod_{i=0}^{p-1} \left[ 
\dim S^{\mu(i)}\right]_{\!p}. 
\end{equation}
Clearly the right-hand side of \eqref{eq:heightineq} is integral, as we required.
\end{proof}


Let $\gamma$ be a $p$-core, let $w\in \mathbf{N}$ and
let $n = |\gamma| + wp$. Let $B$ be the
$p$-block of $S_n$ with $p$-core~$\gamma$ and weight~$w$.
 Applying 
Corollary~\ref{cor:firstRow} to the initial partition $\gamma + wp$ 
defined in \eqref{eq:gamma},
we see that 
there is a vertex of
$S^{\gamma + wp}$ containing a Sylow $p$-subgroup of $S_{wp}$.
Hence~$B$ has a defect group $D$ that contains a Sylow $p$-subgroup
of $S_{wp}$. To complete the proof of Theorem \ref{thm:SymDefectGroups} 
we must show that~$D$ is no larger. 



When defined over fields of characteristic zero, Specht modules
afford the irreducible characters of $S_n$. 
It therefore follows from Theorem~\ref{thm:dimrestriction} that
if~$[n!]_p = p^a$, $p^d$ is the order of the defect group~$D$, 
and $\mu$ is a partition with $p$-core $\gamma$ and weight $w$, then
\begin{equation}
\label{eq:height} 
[\dim S^{\mu}]_p \ge p^{a-d}.
\end{equation}
Moreover, equality holds 
in \eqref{eq:height}
for at least one such partition $\mu$.
But by Lemma~\ref{lemma:weightupperbound}, if~$p^b = [(wp)!]_p$ then
\[ [\dim S^{\mu}]_p \ge p^{a-b} \]
with equality when $\mu = \gamma + wp$. Thus the minimum
in~\eqref{eq:height} occurs when \hbox{$\mu = \gamma + wp$},
and~$p^d = p^b$. Hence $D$ has the same order as a Sylow
subgroup of $S_{wp}$. This completes the proof.






\subsection{On Brauer's Height Zero Conjecture}

Let $G$ be a finite group and let $\chi$ be an ordinary character
of $G$ lying in a $p$-block with defect group of order $p^d$. 
Let $p^a$ be the highest power of $p$ dividing $|G|$. If 
\[ [\chi(1)]_p = p^{a-d+h} \]
then we say that $h$ is the \emph{height} of $\chi$. 
(It follows from Theorem~\ref{thm:dimrestriction} that $h \in \mathbf{N}_0$.)
R.~Brauer made the following
conjecture on character heights.

\begin{conjecture}[Brauer]Every ordinary irreducible character
in a block~$B$ of a finite group has height zero if and only
if  $B$ has an abelian defect group.
\end{conjecture}

Proposition~3.8 in Olsson's paper \cite{OlssonHeights} on character heights and the McKay
Conjecture gives a proof of Brauer's height-zero conjecture
for symmetric groups. It is worth noting
that equation \eqref{eq:heightineq} can be used to give a short alternative
proof.
If $\lambda$ 
is a partition with $p$-quotient $(\mu(0),\ldots,\mu(p-1))$,
where $\mu(i)$ is a partition of $c_i$, then by~\eqref{eq:heightineq},
\[ p^h = 
{\left[ \binom{w}{c_0, 
c_1, \ldots, c_{p-1}} \right]_{\!p}}\, \prod_{i=0}^{p-1} \left[ 
\dim S^{\mu(i)}\right]_{\!p} 
\]
where $h$ is the height of the ordinary character of $S^\lambda$.

If $w < p$ then $c_i < p$ for each $i$, and so each $S^{\mu(i)}$ 
has dimension coprime to $p$. It follows that in blocks of the symmetric group
of abelian defect, every ordinary irreducible character
has height $0$.

If $w \ge p$ then we may 
choose the $c_i$ so that $c_0 + c_1 + \cdots + c_{p-1} = w$ and
the multinomial coefficient is divisible by~$p$,
and then set $\mu(i) = (c_i)$.
Hence in a block of non-abelian defect, there is an ordinary
irreducible character of non-zero height.

\section{The Brauer correspondence for the symmetric group} 
We now 
use the Brauer homomorphism to
determine the Green correspondents of  Specht modules labelled
by initial partitions.
As a corollary of this result (see Corollary~\ref{cor:corr}), we 
get a complete description of how blocks
of symmetric groups relate, under the Brauer correspondence,
to blocks of their local subgroups.

Throughout this section,
let $F$ be a field of characteristic $p$, let $\gamma$ be a
$p$-core and let $w \in \mathbf{N}$. Let $m = |\gamma|$ and
let $n = m + wp$.
The following lemma gives
some of the convenient properties of the Specht modules
$S^{\gamma + wp}$ and $S^\gamma$.

\begin{lemma}\label{lemma:simple}
The Specht module $S^{\gamma + wp}$ is a simple 
$FS_n$-module with trivial source. 
The Specht module $S^{\gamma}$ is a simple projective
$FS_m$-module.
\end{lemma}

\begin{proof}
Because $\gamma + wp$ is the greatest partition labelling
a Specht module in the block with core $\gamma$ and weight $w$, it
follows from \cite[Theorem~12.1]{James} that
$S^{\gamma + wp}$ is simple. For the same
reason, it follows from \cite[Theorem~3.1]{JamesYoung} that
$S^{\gamma + wp}$  is equal to the Young module $Y^{\gamma + wp}$.
Hence $S^{\gamma + wp}$ 
is a direct summand of the permutation module $M^{\gamma + wp}$,
and so it has trivial source. For the second part,
observe that since $\gamma$ is a $p$-core,~$S^\gamma$ is an
indecomposable module lying
in a block of $S_m$ of defect zero. Hence~$S^\gamma$ is simple and projective.
\end{proof}


Our next proposition gives some useful information about
the Brauer quotients
of Specht modules labelled by initial partitions.
In it, we say
that a subgroup $G$ of $S_n$ has 
\emph{support of size $k$} if 
exactly $k$ of the elements
of $\{1,2,\ldots,n\}$ are moved by some permutation in~$G$.

\begin{proposition}\label{prop:Brauercorr}
If $Q$ is a $p$-subgroup
of $S_{wp}$ with support of size $rp$ then
\[ N_{S_n}(Q)/Q \cong S_{n - rp} \times N_{S_{rp}}(Q)/Q \]
and $S^{\gamma + wp}(Q)$ contains a $FN_{S_n}(Q)/Q$-submodule
isomorphic, under this identification, to
\[ S^{\gamma + (w-r)p} \otimes F .\]
If $P$ is a Sylow $p$-subgroup of $S_{wp}$ then
$S^{\gamma + wp}(P) \cong S^{\gamma} \otimes F$.
\end{proposition}

\begin{proof}
As in \S 5, let $t$ be the greatest tableau of shape $\gamma + wp$
in the dominance order. 
By replacing
$Q$ with one of its conjugates if necessary, we can assume
that $Q$ is contained in the symmetric group on the set
\[ Y = \{\gamma_1 + (w-r)p + 1, \ldots,  \gamma_1 + wp \} . \]
Let $X = \{ 1,2,\ldots,n\} \backslash Y$.
Clearly we have
$N_{S_n}(Q) = S_X \times N_{S_Y}(Q)$, which implies
the first assertion in the theorem.


Let $U$ be the kernel of the Brauer homomorphism
from $(S^{\gamma + wp})^Q$ to $S^{\gamma + wp}(Q)$.
Since $Q$ is contained in the subgroup $H(t)$
defined in Theorem~\ref{thm:main}, it follows from
Lemma~\ref{lemma:kernel} that 
$e_t + U$ is a non-zero
element of $S^{\gamma + wp}(Q)$. Let~$W$ 
be the submodule of $S^{\gamma + wp}(Q)$ generated by $e_t + U$.
By Lemma~\ref{lemma:fixed},~$e_t$ is fixed by every permutation in $S_Y$, 
so 
$N_{S_Y}(Q)/Q$ acts trivially on $W$.

Let $s$ be the greatest tableau of shape
$\gamma + (w-r)p$ in the dominance order.
Restricting to the action of $S_X \cong S_{n-rp}$, we see that there
is a surjective map of $FS_{n-rp}$-modules,
\[ S^{\gamma + (w-r)p} \rightarrow W \] 
defined by extending $e_s \mapsto e_t + U$. This map is non-zero
since $e_t \not\in U$. Moreover, by 
Lemma~\ref{lemma:simple}, $S^{\gamma + (w-r)p}$
is a simple $FS_{n-rp}$-module. Hence 
\[ W \cong S^{\gamma + (w-r)p} \otimes F \]
as a module for $FN_{S_n}(Q)/Q$. 

It only remains to show that if 
$P$ is a Sylow
$p$-subgroup of $S_Y \cong S_{wp}$, 
then $S^{\gamma + wp}(P) \cong S^{\gamma} \otimes F$.
By  Lemma~\ref{lemma:simple}, $S^\gamma$ is a projective
$FS_m$-module, and since $N_{S_Y}(P)/P$ has order coprime to $p$,
the trivial $FN_{S_Y}(P)/P$-module is projective.
Hence $S^\gamma \otimes F$ 
is a projective as a module for 
$F(S_X \times N_{S_Y}(P))$. By the previous paragraph,
it splits off as a direct summand of $S^{\gamma + wp}(P)$.
Since~$S^{\gamma + wp}$ has trivial source,
it follows from Theorem~\ref{thm:Brauerpperm} that
$S^{\gamma  + wp}(P)$ is indecomposable (and projective)
as an $FN_{S_n}(P)/P$-module. Therefore $S^{\gamma + wp}(P) = W$,
as claimed.
\end{proof}

By Theorem~\ref{thm:Brauerpperm}(ii), 
the Green correspondent of a  trivial source module
is equal to its Brauer quotient.
Proposition~\ref{prop:Brauercorr} therefore
implies the following theorem, which gives
a complete description of
the local properties of Specht modules labelled by initial partitions.
It should be noted that 
in \cite{Valero} the author proves Theorem~\ref{thm:Greencorr}
(by an explicit calculation) in the case when $p=2$ and~$\gamma$ is 
a $2$-core. 

\begin{theorem}\label{thm:Greencorr}
Let $P$ be a Sylow $p$-subgroup of $S_{wp}$.
The Specht module $S^{\gamma + wp}$ has $P$ as
one of its vertices, and its source is the trivial $FP$-module.
The Green correspondent of $S^{\gamma + wp}$
is the $FN_{S_{n}}(P)$-module $S^\gamma \otimes F$.
\end{theorem}


We end by using Proposition~\ref{prop:Brauercorr} to
describe 
how blocks of symmetric groups behave under the Brauer
correspondence. (For the original proof by
M.~Brou{\'e}, see \cite{BroueBlocks}.) Our 
definition of the Brauer correspondence
is taken from Alperin \cite[\S 14]{Alperin};
thus if $H$ is a subgroup of a finite group $G$ and $b$ is a block
of $H$, then we say that 
$b$ \emph{corresponds to the block}~$B$ of $G$,
and write $b^G = B$, if $b$, considered
as a $F(H \times H)$-module, is a 
summand of the restriction of $B$ to $H \times H$, and~$B$
is the unique block of~$G$ with this property.
 
We shall need the following lemma, which 
generalises a well-known result  about the Green correspondence
(see, for example \cite[\S14, Corollary 4]{Alperin}) to
Brauer quotients. 

\begin{lemma}\label{lemma:genAlperin}
Let $G$ be a finite group and
let $M$ be an indecomposable $FG$-module with vertex $P$ and trivial source.
Let $Q$ be a subgroup of $P$.
Suppose that  $M$ lies in the block~$B$ of~$G$. If~$M(Q)$, considered
as an $FN_G(Q)$-module,
has a summand in the block~$b$ of~$N_G(Q)$, 
then 
$b^G$ is defined and $b^G= B$.
\end{lemma}

\begin{proof}
By \cite[Exercise~27.4]{Thevenaz},
when considered as an $FN_G(Q)$-module, $M(Q)$ has a summand whose 
vertex contains~$Q$. Hence there is 
some defect group~$D$ of~$b$ which contains $Q$. Therefore $N_G(Q) \supseteq
C_G(Q) \supseteq C_G(D)$, and so by part~3 of Lemma~1 on 
page~101 of~\cite{Alperin},~$b^G$ is defined. 

Again by \cite[Exercise~27.4]{Thevenaz}, 
$M(Q)$ is a direct summand of $M\res_{N_G(Q)}$.
Hence $M(Q)$ is not killed by $B \res_{N_G(Q) \times N_G(Q)}$, and 
so 
\[ B \res_{N_G(Q) \times N_G(Q)}b \not= 0. \] 
Hence $b^G = B$, as required.
\end{proof}

\begin{corollary}\label{cor:corr}
Let $Q$ be a $p$-subgroup of $S_{wp}$ with support of size $rp$.
In the Brauer correspondence 
between blocks of $S_n$ and blocks of 
$N_{S_n}(Q) \cong S_{n-rp} \times N_{S_{rp}}(Q)$,  
the $p$-block of $S_n$ with core $\gamma$
and weight $w$ corresponds to 
\[ b \times b_0(N_{S_{rp}}(Q)) \]
where $b$ is the $p$-block of $S_{n-rp}$ with core $\gamma$ and weight $w-r$,
and $b_0(N_{S_{rp}}(Q))$ is the principal block of $N_{S_{rp}}(Q)$.
\end{corollary}

\begin{proof}
It follows from Proposition~\ref{prop:Brauercorr} that 
$S^{\gamma + wp}(Q)$ has a summand in the block 
$b \times b_0(N_{S_{wp}}(Q))$. Now apply
Lemma~\ref{lemma:genAlperin}.
\end{proof}

In his earlier proof of Corollary~\ref{cor:corr},
Brou{\'e} notes that the group $N_{S_{rp}}(Q)$ has a unique $p$-block 
(see \cite[Lemma~2.6]{BroueBlocks}).
The correspondence described in
Corollary~\ref{cor:corr} is therefore bijective. 


\section{The vertices of $S^{(n-2,2)}$}

We end by using Theorem~\ref{thm:main} to find 
the vertices of $S^{(n-2,2)}$ over fields
of odd characteristic. We shall need Corollary~1 of \cite{Green59},
which states that if $G$ is a finite group,~$F$ is a field
of characteristic $p$, and $M$ is an indecomposable $FG$-module
of dimensional coprime to $p$, then $M$ has a Sylow $p$-subgroup of $G$
as one of its vertices. Our result is as follows.

\begin{theorem}
Let $n \ge 4$.
When defined over a
field of odd characteristic~$p$,
the Specht module $S^{(n-2,2)}$ is indecomposable, and its vertex
is equal to the defect group of the $p$-block in which it lies.
\end{theorem}

\begin{proof}
That $S^{(n-2,2)}$ is indecomposable follows from Theorem~\ref{thm:indec}.
The dimension of $S^{(n-2,2)}$ is $n(n-3)/2$. Hence if neither
$n$ nor $n-3$ is divisible by $p$, then, by the result of Green
just mentioned, $S^{(n-2,2)}$ has a Sylow $p$-subgroup of $S_n$ as 
its vertex. This is also a defect group of its $p$-block. 

If $p$ divides $n$ and $p > 3$ then $(n-2,2)$ has $p$-core $(p-2,2)$
and weight $(n-1)/p$. Hence the defect group of its block is a Sylow
$p$-subgroup of $S_{n-p}$. It
follows from Theorem~\ref{thm:main}
that
$S^{(n-2,2)}$ has a vertex containing a Sylow $p$-subgroup of $S_{n-p}$.
Therefore the vertex of $S^{(n-2,2)}$ agrees with its defect group.
If $p = 3$ then
$(n-2,2)$ lies in the block with $p$-core $(4,2)$ and weight $(n-6)/3$.
Theorem~\ref{thm:main} implies that there is a vertex containing
a Sylow $3$-subgroup of $S_{n-6}$, as required.

The only remaining case is when $p$ divides $n-3$ and $p > 3$. 
In this case $(n-2,2)$ has $p$-core $(p+1,2)$ and weight $(n-3)/p\: - 1$,
and the result follows from Theorem~\ref{thm:main} as before.
\end{proof}

\section*{Acknowledgement}

This paper is based on \S 3.5 of the author's D.~Phil thesis,
written under the supervision of Karin Erdmann. I gratefully
acknowledge her advice and support.



\begin{thebibliography}{10}

\bibitem{Alperin}
{\sc Alperin, J.~L.}
\newblock {\em Local representation theory}, vol.~11 of {\em Cambridge studies
  in advanced mathematics}.
\newblock Cambridge University Press, 1986.

   

\bibitem{AlperinWeights}
{\sc Alperin, J.~L.}
\newblock Weights for finite groups.
\newblock In \emph{The {A}rcata {C}onference on {R}epresentations of {F}inite
{G}roups ({A}rcata, {C}alif., 1986)},
{\em Proc. Sympos. Pure Math. 47.} 
\newblock Amer. Math. Soc. (1987), 369--379.

\bibitem{Broueperm}
{\sc Brou{\'{e}}, M.}
\newblock On {S}cott modules and {$p$}-permutation modules: an approach through
  the {B}rauer homomorphism.
\newblock {\em Proc.~Amer.~Math.~Soc. 93} (1985), 401--408.

\bibitem{BroueBlocks}
{\sc Brou{\'{e}}, M.}
\newblock Les {$l$}-blocs des groupes {$GL(n,q)$} et {$U(n,q^{2})$} et leurs
  structures locales.
\newblock {\em Ast{\'{e}}risque 133--134\/} (1986), 159--188.

\bibitem{Cameron}
{\sc Cameron, P.~J.}
\newblock {\em Combinatorics: Topics, Techniques, Algorithms}.
\newblock Cambridge University Press, 1994.

\bibitem{DanzErdmann}
{\sc Danz, S., and Erdmann, K.}
\newblock The vertices of a class of {S}pecht modules and simple modules for
  symmetric groups in characteristic 2.
\newblock {\em Preprint\/} (May 2009). Available from
\url{http://www.minet.uni-jena.de/algebra/preprints/Danz/twopartchar2.pdf}

\bibitem{DKZ}
{\sc Danz, S., K{\"u}lshammer, B. and Zimmermann, R.}
\newblock On vertices of simple modules for symmetric groups of small degrees.   
\newblock {\em J. Algebra 320} (2008), 680--707.   

\bibitem{Dornhoff}
{\sc Dornhoff, L.}
\newblock {\em Group representation theory, Part II}.
\newblock Dekker, 1972.

\bibitem{Green59}
{\sc Green, J.~A.}
\newblock On the indecomposable representations of a finite group.
\newblock {\em Math. Zeitschrift 70\/} (1958/59), 430--445.

\bibitem{JK}
{\sc James, G., and Kerber, A.}
\newblock {\em The representation theory of the symmetric group}, vol.~16 of
  {\em Encyclopedia of Mathematics and its Applications}.
\newblock Addison-Wesley Publishing Co., Reading, Mass., 1981.

\bibitem{JamesCounters}
{\sc James, G.~D.}
\newblock Some counterexamples in the theory of {S}pecht modules.
\newblock {\em J. Algebra 46} (1977), 457--461.

\bibitem{James}
{\sc James, G.~D.}
\newblock {\em The representation theory of the symmetric groups}, vol.~682 of
  {\em Lecture Notes in Mathematics}.
\newblock Springer, Berlin, 1978.

\bibitem{JamesYoung}
{\sc James, G.~D.}
\newblock Trivial source modules for symmetric groups
\newblock {\em Arch.~Math. (Basle) 41} (1983), 294--300.

\bibitem{MZ}
{\sc M{\"u}ller J., and Zimmermann, R.}
\newblock Green vertices and sources of simple modules of the symmetric group.
\newblock {\em Arch. Math. (Basle) 89} (2007), 87--108.

\bibitem{Murphy}
{\sc Murphy, G.~E.}
\newblock A new construction of Young's seminormal representation of the
symmetric group.
\newblock {\em J. Algebra 69} (1981), 287--297.

\bibitem{MurphyDec}
{\sc Murphy, G.~M.}
\newblock On decomposability of some {S}pecht modules for symmetric groups.
\newblock {\em J. Algebra 66} (1980), 156--168.

\bibitem{MurphyPeel}
{\sc Murphy, G.~M., and Peel, M.~H.}
\newblock Vertices of {S}pecht modules.
\newblock {\em J. Algebra 86} (1984), 85--97.

\bibitem{OlssonHeights}
{\sc Olsson, J.~B.}
\newblock Mc{K}ay numbers and heights of characters.
\newblock {\em Math. Scand. 38}  (1976), 24--42.

\bibitem{Peel75}
{\sc Peel, M.~H.}
\newblock Specht modules and symmetric groups.
\newblock {\em J. Algebra 36\/} (1975), 88--97.

\bibitem{Specht1935}
{\sc Specht, W.}
\newblock Die irreduziblen {D}arstellugen der symmetrischen {G}ruppe.
\newblock {\em Math. Zeitschrift 39} (1935), 696--711.

\bibitem{Thevenaz}
{\sc Th{\'e}venaz, J.}
\newblock {\em G-Algebras and modular representation theory}.
\newblock Oxford mathematical monographs. Oxford University Press, 1995.

\bibitem{Valero}
{\sc Valero-Elizondo, L.}
\newblock Triangular partitions with augmented first rows and weights for the
  symmetric groups in characteristic two.
\newblock {\em Bol. Soc. Mat. Mexicana (3) 8} (2002), 117--126.

\bibitem{WildonCycSpechtVertices}
{\sc Wildon, M.}
\newblock Two theorems on the vertices of {S}pecht modules.
\newblock {\em Arch. Math. (Basel) 81} (2003), 505--511.

\bibitem{WildonThesis}
{\sc Wildon, M.}
\newblock {\em Modular representations of symmetric groups}.
\newblock D.~Phil Thesis, Oxford University, 2004.

\end{thebibliography}

\def\cprime{$'$} \def\Dbar{\leavevmode\lower.6ex\hbox to 0pt{\hskip-.23ex
  \accent"16\hss}D}


\end{document}